\documentclass[12pt]{article}

\usepackage{amsmath}

\usepackage{amssymb}

\usepackage{exscale}

\usepackage{amsthm}

\usepackage{epsfig}

\usepackage{verbatim}

\theoremstyle{plain}

\newtheorem{theorem}{Theorem}

\begin{document}

\title{Laplace spectra on open and compact\\ Zeeman manifolds}

\author{Zolt\'an I. Szab\'o}

\date{ }

\maketitle

\begin{abstract}
By a recent observation, the Laplacians on the
Riemannian manifolds the author used for isospectrality 
constructions are nothing but the 
Zeeman-Hamilton operators of free 
charged particles. These manifolds can be considered as prototypes 
of the so called Zeeman manifolds. This observation allows to 
develop a spectral theory
both on open Z-manifolds and their compact submanifolds.

The theory on open manifolds leads to a new nonperturbative approach
to the infinities of QED.  This idea exploits that 
the quantum Hilbert space, 
$\mathcal{H}$, 
decomposes into subspaces (Zeeman zones) which are 
invariant under the actions
both of this Zeeman-Laplace operator and the natural Heisenberg group 
representation. Thus a well defined particle 
theory and zonal geometry can be developed on each zone separately. 
The most surprising result
is that quantities divergent
on the global setting are finite on the
zonal setting. Even the zonal Feynman integral is well defined.
The results include explicit computations of objects
such as the zonal spectra, the waves defining the zonal point-spreads, 
the zonal Wiener-Kac resp. Dirac-Feynman flows, and the corresponding
partition functions. 

The observation adds new view-point also to the problem of finding 
intertwining operators by which isospectral pairs of metrics
with different local geometries 
on compact submanifolds can be constructed. Among the
examples the author constructed the most surprising are the  
isospectrality families containing both 
homogeneous and locally inhomogeneous metrics. The 
observation provides even quantum physical interpretation to the
isospectrality.
\end{abstract}

\section{Zeeman manifolds}

{\bf Zeeman-Hamilton operators.} 
The classical Zeeman operator of a charged particle is 
\begin{equation}
\label{Zee_int}
H_Z=-{\hbar^2\over 2\mu}\Delta_{(x,y)} -
{\hbar eB\over 2\mu c\mathbf i}
D_z\bullet
+{e^2B^2\over 8\mu c^2}(x^2+y^2) +eV ,
\end{equation}
where $V$ is the Coulomb potential originated from the nucleus 
(for free particles $V=0$ holds).
At the time it was introduced, the new feature of this operator was 
the orbital angular momentum operator 
$D_z\bullet =x\partial_y-y\partial_x$, which was the forerunner for
an adequate spin-concept. 
This operator is the result of a long agonizing creative 
process \cite{to},
which was used for explaining 
the Zeeman effect. (Note that the $D_z\bullet$
commutes with the rest part, $\bold O$, of the complete operator.
Thus the spectrum appears on common eigenfunctions, resulting that the
$D_z\bullet$ splits the spectral lines of $\bold O$ (Zeeman effect).) 
Pauli, who added a spin angular momentum operator to the orbital one,
developed the non-relativistic spin-concept. The relativistic concept
due to Dirac. 
Actually, the $H_Z$ is the Hamilton operator of an 
electron orbiting about the origin of the 
$(x,y)$-plane in a constant magnetic field $\mathbf K=B\partial_z$.
It had been established by means of the Maxwell equations.
\medskip

{\bf Mathematical modeling: Zeeman manifolds.}
An interesting feature of free Zeeman operators, $H_Z$, is that they 
can be identified 
with the Laplace operators of certain Riemannian manifolds, namely,
with the Laplacians on two step nilpotent Lie groups endowed 
with the natural left invariant metrics. The details are as follows.

A 2-step nilpotent metric Lie group is defined on the
product $\mathbf v\oplus\mathbf z$ of Euclidean spaces, where the
components, 
$\mathbf v=\mathbb R^k$ and
$\mathbf z=\mathbb R^l$, 
are called X- and Z-space respectively. The
main object defining the Lie algebra is the linear space $J_{\mathbf z}$
of skew endomorphisms $J_Z$, where $Z\in\mathbf z$, acting on the
X-space. The metric, $g$, is the left invariant extension of the
natural Euclidean metric on the Lie algebra. 

Particular 2-step nilpotent Lie groups are the
Heisenberg-type Lie groups which are defined by endomorphism spaces
satisfying the Clifford condition $J^2_Z=-|Z|^2id$. I. e.,
these metric groups are attached to Clifford modules.
In this case the X-space decomposes into the product 
$\bold v=(\mathbb R^{r(l)})^{a+b}=\mathbb R^{r(l)a}\times
\mathbb R^{r(l)b}$.
Endomorphisms $J_Z$ 
are defined by endomorphisms $j_Z$
acting on the smaller space $\mathbb R^{r(l)}$. Namely, the $J_Z$
acts on  
$\mathbb R^{r(l)a}$ resp.
$\mathbb R^{r(l)b}$ as $j_Z\times\dots\times j_Z$ 
resp. $-j_Z\times\dots\times -j_Z$. 
The H-type groups are denoted by 
$H^{(a,b)}_l$, indicating the above decomposition. A point on the
nilpotent group is represented by $(X,Z)$.
Each group, $(N,g)$, extends into a solvable group $(SN,g_s)$,
where a point is represented by $(t,X,Z)$. 
The solvable extensions of H-type groups are then $SH^{(a,b)}_l$.  

The Laplacians on H-type groups are of the form
\begin{equation}
\label{Delta}
\Delta=\Delta_X+(1+\frac 1{4}|X|^2)\Delta_Z
+\sum_{\alpha =1}^r\partial_\alpha D_\alpha \bullet,
\end{equation}
where $D_\alpha\bullet$ denotes directional derivatives along
the X-fields 
$J_\alpha (X)=J_{Z_\alpha}(X)$ and endomorphisms $J_\alpha$ are 
defined by an orthonormal
basis $\{ Z_\alpha\}$ of the Z-space. 

This operator is not the
Zeeman operator yet. The Zeeman operator appears on 
center periodic H-type groups  
introduced by factorizations, $\Gamma_\gamma\backslash H$, 
with Z-lattices $\Gamma_\gamma =\{Z_\gamma\}$ of the Z-space. 
In fact, in this case the $L^2$ function
space is the direct sum of function
spaces $W_\gamma$ spanned by functions of the form
$\Psi_\gamma (X,Z)
=\psi (X)e^{2\pi\mathbf i\langle Z_\gamma ,Z\rangle}.
$
Each $W_\gamma$ is invariant under the action of the Laplacian, i. e., 
$\Delta \Psi_{\gamma }(X,Z)=\Box_{\gamma}\psi (X)
e^{2\pi\mathbf i\langle Z_\gamma ,Z\rangle}$,  
where operator $\Box_{\gamma }$, acting on 
$L^2(\mathbf v)$, is of the form
\begin{equation}
\label{Box}
\Box_{\gamma }
=\Delta_X + 2\pi\mathbf i D_{\gamma }\bullet -4\pi^2
|Z_\gamma |^2(1 + 
\frac 1 {4} |X|^2).
\end{equation}
  
Notice that (\ref{Zee_int}) is nothing but (\ref{Box}) 
on the 3D-Heisenberg
group. On a $(k+1)$-dimensional Heisenberg group, defined by a complex
structure $J$ acting on the even dimensional Euclidean space  
$\mathbf v=\mathbf R^k$, the Laplacian (\ref{Box}) appears in the form 
\begin{equation}
\label{Box2}
\Box_{\lambda}=\Delta_X +2 \mathbf i 
D_{\lambda }\bullet -
\lambda^2|X|^2-4\lambda^2.
\end{equation}
Number $k/2$ is interpreted as the number of particles. The single
complex structure $J$ is interpreted such 
that these particles
are rotating in the same plane defined by the same 
constant magnetic field $B$. In the 2D-version the particles
belonging to $J$ and $-J$ are called antiparticles (electrons and
positrons). The single parameter $\lambda$ above is interpreted
that the particles are identical up-to the sign of the charge they
are loaded.

Operators (\ref{Zee_int}) and (\ref{Box2}) are identified by 
$H_Z=-(1/2)\Box_\lambda$ and by the particular choice 
$\mu =\hbar =1, \lambda=eB/2c$ of the constants.
Operator (\ref{Box}) contains also the
constant $-4\pi^2|Z_\gamma |^2=4\lambda^2$, which is proportional
to $\hbar^2$ on the microscopic level, thus, it is usually neglected
in quantum physics. Also note that
the particles described by these Hamiltonians are free ($V=0$).

The Zeeman operator appears as Laplacian on center periodic
2-step nilpotent Lie groups in a more complex form.
These models represent $k/2$ number of charged particles, each of them
is orbiting in its own constant magnetic field. The system
can be in crystal states represented by the endomorphisms
$J_{\gamma}$. The Hamilton operators belonging to these crystal states
are $-{1\over 2}\Box_\gamma$.  This model matches
Dirac's famous multi-time theory. 

The Riemannian manifolds introduced so far are prototypes of a general
Zeeman manifold concept. This general concept
is beyond the scope of this talk and will be developed in a 
subsequent paper. This exposition proceeds with considering
2-step nilpotent Lie groups.
\medskip

\section{Isospectrality constructions}

The isospectrality constructions are performed on H-type
groups, 
$H^{(a,b)}_l$, 
and on their solvable extensions, 
$SH^{(a,b)}_l$, 
both with
periodic and non-periodic centers. For fixed $a+b$ and $l$, these
groups are defined on the same $(X,Z)$- resp. $(t,X,Z)$-space.
There is established in many different ways that the local
geometries regarding the metrics $g_l^{(a,b)}$ of a family 
are completely different (cf. the striking examples). Yet, 
certain submanifolds of
these common underlying manifolds have the same Laplace spectra
regarding all the metrics of the family. The endomorphism space
$J_{l}^{(a^\prime,b^\prime)}$
can be constructed by an other endomorphism space,
$J_{l}^{(a,b)}$, 
of the family such that the irreducible endomorphism spaces,
$j_{l}$ or 
$-j_{l}$, 
by which $J_l$ is expressed are multiplied by $-1$ on some irreducible
subspaces, $\mathbb R^{r(l)}$. In retrospect to the physical
interpretation, this means that some of the particles are exchanged
by their antiparticles. Although this exchange drastically changes
the local geometry, the spectrum on the considered domains remains
the same.

{\bf Constructing the ball$\times$torus- and 
sphere$\times$torus-examples.} These examples are constructed for a
family,
$\Gamma\backslash H^{(a,b)}_l$, 
of Z-periodic manifolds
defined for fixed $l$ and $a+b$. The 
submanifolds considered are torus bundles over a ball (resp. sphere) 
around the origin of the X-space.

To establish an intertwining operator, one should consider,
for each invariant space $W_\alpha$ constructed above, 
an orthogonal transformation on the X-space which conjugate 
$J_{\alpha}^{(a,b)}$ to
$J_{\alpha}^{(a^\prime,b^\prime)}$. This transformation clearly
intertwines $\Box_\alpha$ with $\Box^\prime_\alpha$ such that it
keeps also the boundary conditions. It induces an appropriate
intertwining operator also on the boundary manifolds.
The striking examples appear on the quaternionic families 
$H^{(a,b)}_3$. In this case the 
sphere$\times$torus-type boundary manifolds in 
$\Gamma\backslash H^{(a+b,0)}_3$
are homogeneous while the others in the isospectrality family are
locally inhomogeneous. 

In retrospect, this intertwining operator was constructed in two steps:
First, the $L^2$-function space has been decomposed into invariant
subspaces, then, for each invariant subspace an operator intertwining
both the Laplacians and the boundary conditions has been found. This
operator is derived from a point transformation. For distinct
invariant subspaces also these point transformations are distinct. 
Note that for each invariant
space only a single endomorphism $J_\gamma$ is involved to the
constructions. 
 
{\bf Constructing the ball- and sphere-type examples.} These examples
were originally constructed in \cite{sz2,sz3}. It turned out just
recently that the intertwining operator was incorrectly defined in
these articles. The corrected intertwining operator has been introduced
in \cite{sz4}. This correction saved all the previous results and,
additionally, it produced new examples defined on 
sphere$\times$ball- and
sphere$\times$sphere-type manifolds. Among them also new striking
examples have been found. Next this corrected intertwining operator
is constructed.
  
The ball-type domains
are, by definition, diffeomorphic to Euclidean balls
such that the sphere-type boundary manifolds
are level sets described by equations of the form 
$\varphi (|X|,|Z|)=0$ resp. $\varphi (|X|,|Z|,t)=0$. 
These domains are invariant under the action of the orthogonal
group $\mathbb O(\mathbb R^k)\times\mathbb O(\mathbb R^l)$ (domains of 
$(X,Z)$-revolution), which can be visualized such that
there is a ball, centered at the origin in the X-space, over the 
points of which
there are Z-balls of radius $R_Z(|X|)$ around the origin
in the Z-space. Then the boundary is a level set defined by
$\varphi (|X|,|Z|)=|Z|-R_Z(|X|)=0$. 

Note that radius $R_Z(|X|)$
is constant along a sphere defined by a constant radius $R_X=|X|$
in the X-space and the
ball bundle defined by the Z-balls 
over this X-sphere is trivial. These are
the so called {\it sphere$\times$ball-type manifolds} whose boundaries
are {\it sphere$\times$sphere-type manifolds}. The new examples,
not discussed in the earlier papers, are constructed on these
domains and surfaces.

This visualization
can be started out by a Z-ball in the Z-space over the points of
which there are X-balls of radius $R_X(|Z|)$ given.
Then the boundary is defined by $|X|-R_X(|Z|)=0$. This paper
proceeds with the first approach.

In the solvable
case one should consider $(Z,t)$-balls and $(Z,t)$-spheres around 
the origin $(0_Z,1)$. The base manifold is the same X-sphere as before.
Note that a Z-ball $B_{R_Z}(0_Z)$ (resp. Z-sphere $S_{R_Z}(0_Z)$)
uniquely extends into a geodesic ball (resp. sphere) of the  
(hyperbolic) $(Z,t)$-space. 
A sphere$\times$ball-type domain can be regarded as a level set
in a ball-type domain such that the Z-balls (or $(Z,t)-balls)$ of the
ball type domain are considered over a sphere $S_{R_X}$ in the X-space. 
Similarly, the sphere$\times$sphere-type manifolds can be regarded
as level sets in sphere-type manifolds.

The isospectrality will be
investigated in details for the discrete families, 
$H^{(a,b)}_l$, 
of Heisenberg type groups 
defined by the same $l$ and $a+b$. The Laplacian to be considered now
is described in (\ref{Delta}). Comparing with the Zeeman operator
(\ref{Box}), this operator involves all the endomorphisms, making
the constructions rather difficult. The Laplacians of the members in
a family differ from each other just by the last term, $\bold M$,
which is a compound orbital angular momentum operator. 
There is explained
in \cite{sz5} that the compound orbital angular momentum operator 
in Dirac's famous multi-time
model, which complies with relativity by endowing each particle in the
system with self-time, appears in this form. The spectral
investigation both of $\bold M$ and $\Delta$ is completely missing
both in physics and mathematics. This explains some extend the 
difficulties one should face in investigating this Laplacian.
Note that this most intriguing operator, $\bold M$,  
commutes with both operators in the rest part of (\ref{Delta}). 

{\bf Constructing the intertwining operator.} 
In the {\it first step} we 
explicitly describe the {\it eigenfunctions both of $\bold M$ and
$\Delta$ with no boundary conditions}. Since the $\bold M$
commutes with the rest part $\bold O$ of $\Delta$, the 
eigenfunctions can be
sought such that they are eigenfunctions 
both of $\bold M$ and $\bold O$.

In the very first step we look for the eigenfunctions of a single 
angular momentum operator $\bold D_V\bullet$, defined for a
Z-vector $V$.  For a fixed X-vector $Q$
and the unit Z-vector $V_u={1\over |V|}V$, consider the X-function
$\Theta_{Q}(X,V_u)=\langle Q+\bold iJ_{V_u}(Q),X\rangle$ and its 
conjugate $\overline{\Theta}_{Q}(X,V_u)$. For vector $V=|V|V_u$,
these functions are 
eigenfunctions of $D_{V}\bullet$ with eigenvalue $-|V|\bold i$ resp.
$|V|\bold i$. The higher order eigenfunctions are of the form
$\Theta_{Q}^p\overline\Theta^q_{Q}$ 
with eigenvalue $(q-p)|V|\bold i$.
 
In order to find the eigenfunctions of the compound operator $\bold M$,
consider a sphere $S_{R_Z}$ of radius $R_Z$
around the origin in the Z-space. For an appropriate function
$\phi (|X|,V)$, depending on $|X|$ and $V\in S_{R_Z}$,
define 
\begin{equation}
\label{FourR_Z}
\mathcal F_{QpqR_Z}(\phi )(X,Z)
=\oint_{S_{R_Z}}e^{\bold i
\langle Z,V\rangle}\phi (|X|,V)
\Theta_{Q}^p(X,V_u)\overline\Theta^q_{Q}(X,V_u)dV.
\end{equation} 
By $\bold M\oint =\oint \bold iD_V\bullet$, 
this function is an eigenfunction
of $\bold M$ with the real eigenvalue $(p-q)R_Z$. 
These functions are eigenfunctions also of $\Delta_Z$ with eigenvalue
$R_Z^2$. Also note 
that these eigenvalues do not change by varying $Q$.

The function space spanned by functions (\ref{FourR_Z}) 
generated by different $\phi$'s is not invariant with respect to 
the action of $\Delta_X$, thus the eigenfunctions of the complete
operator $\Delta$ do not appear in this form. In order to find the
common eigenfunctions,
the homogeneous non-harmonic polynomials 
$
\Theta_{Q}^p\overline\Theta^q_{Q}
$
of the X-variable should be exchanged in the above formula for the 
polynomials
$
\Pi_X (\Theta_{Q}^p\overline\Theta^q_{Q})
$,
defined by projections, $\Pi_X$, onto the space of $(p+q)$-order 
homogeneous harmonic 
polynomials of the X-variable. These projections 
are described in the form 
$\Pi_X =\Delta_X^0+B_1|X|^2\Delta_X +B_2|X|^4\Delta_X^2+\dots$
in \cite{sz3} (cf. (3.14)). 
By this formula easily follows that also the functions 
\begin{equation}
\label{HFourR_Z}
\mathcal {HF}_{QpqR_Z}(\phi )(X,Z)
=\oint_{S_{R_Z}}e^{\bold i
\langle Z,V\rangle}\phi (|X|,V)
\Pi_X(\Theta_{Q}^p(X,V_u)\overline\Theta^q_{Q}(X,V_u))dV
\end{equation}
are eigenfunctions
both of $\bold M$ and $\Delta_Z$ with the same eigenvalues as in
(\ref{FourR_Z}). 

When the complete Laplacian (\ref{Delta}) is
acting on this function, one gets an action which is combination of
X-radial differentiation, $\partial_{|X|}$, 
and multiplications with functions depending 
just on $|X|$. I. e., the action is completely reduced to X-radial
functions. Also this reduced form of the Laplacian is not changing by 
varying $Q$. The eigenfunctions of $\Delta$ can be found by seeking
the eigenfunctions of the reduced operator among the X-radial functions.
Since it is not needed, we do not go into further details 
of this computation here.

The above discussions clearly indicate that functions (\ref{FourR_Z}), 
defined for all $S_{R_Z}$, can be used to define {\it operators
intertwining the Laplacians} for any two Heisenberg type groups
$H^{(a,b)}_l$ and $H^{(a^\prime,b^\prime)}_l$ satisfying
$a+b=a^\prime+b^\prime$. The intertwining operators 
$
\kappa_{Qpq}   
$
can be loosely defined by the correspondence
\begin{equation}
\kappa_{Qpq}:  
\mathcal F_{QpqR_Z}(\phi )(X,Z)\to
\mathcal {F}^{\prime}_{QpqR_Z}(\phi )(X,Z),
\end{equation}
where the associated functions are built up by
$J_\alpha$ resp. $J^\prime_\alpha$ in the very same form. For fixed
$Q$ and $n,m\in \mathbb N$, functions
$
\mathcal F_{QpqR_Z}(\phi )(X,Z)
$ 
defined for $p+q\leq n\, ,\,q-p=m$ and all $R_Z$ span a function space
which is closed with respect to the action of $\Delta$. The
$
\kappa_{Qpq}  
$
maps this function space to the corresponding other one by 
intertwining $\Delta$ and $\Delta^\prime$ term by term. 
The intertwining of $\bold M$ resp. $\Delta_Z$ with $\bold M^\prime$
resp. $\Delta_Z$ is clear by the above considerations. One can show
that also the 
$\Delta_X$ is intertwined with itself. Moreover, the correspondence  
$
\kappa_{Qpq}:  
\mathcal{HF}_{QpqR_Z}(\phi )(X,Z)\to
\mathcal {HF}^{\prime}_{QpqR_Z}(\phi )(X,Z),
$
defines the same operator as before.
Actually, the function spaces of the first type are direct sums
of spaces of the second type. The operators defined by the two 
alternatives agree. The second definition clearly shows the
intertwining of $\Delta_X$ and $\Delta_X^\prime$. 

This was just a loose description of the intertwining operator.
One should also keep in mind that it was considered just on an
invariant subspace defined by a fixed $Q$ (the global intertwining
will be defined, in the end, by splitting the $L^2$-space into
such invariant spaces). 
For defining this reduced intertwining in a precise form, 
notice that (\ref{FourR_Z}) 
describes the Fourier transform of
a Dirac-type distribution concentrated on the sphere $S_{R_Z}$.
In polar coordinates the Z-Fourier transform can be interpreted as  
superposition of 
transforms described in (\ref{FourR_Z}). Thus the map defined
by the Z-Fourier transform
\begin{eqnarray}
\kappa_{Qpq} : &
\mathcal F_{Qpq}(\phi )=
\int_{\mathbb R^l} e^{\bold i\langle Z,V\rangle}
\phi (|X|,V)\Theta_{Q}^p(X,V_u)\overline\Theta^q_{Q}(X,V_u)dV\to \\
& \to \mathcal F^\prime_{Qpq}(\phi )=
\int_{\mathbb R^l} e^{\bold i\langle Z,V\rangle}
\phi (|X|,V)\Theta_{Q}^{\prime p}(X,V_u)
\overline\Theta_Q^{\prime q}(X,V_u)dV.
\nonumber 
\end{eqnarray}
is indeed an intertwining operator (cf. also (\ref{MF})-(\ref{HMF})).
 
The $L^2$-function spaces between which this intertwining is defined
are denoted by 
$
\bold \Phi_{Qpq}
$ and
$\bold \Phi^{\prime}_{Qpq}$ respectively. They are the Z-Fourier
transforms of the $L^2$-Hilbert spaces spanned by functions
of the form 
$\phi\Theta^p_Q\overline\Theta^q_Q$ and
$\phi\Theta^{\prime p}_Q\overline\Theta^{\prime q}_Q$, 
where $\phi(Z)$
is an arbitrary $L^2$-function. 
The corresponding function spaces defined by the Fourier transform
of functions involving 
the X-harmonic polynomials 
$\Pi_X(\Theta_Q^p\overline\Theta_Q^q)$ are denoted by 
$\bold\Xi_{Qpq}$.
They are corresponded to 
$\bold\Xi^{\prime}_{Qpq}$ by the 
$\kappa_{Qpq}: 
\mathcal{HF}_{Qpq}(\phi )
\to \mathcal{HF}^\prime_{Qpq}(\phi )$. (It should not be confusing
that this map is denoted by the same formula regarding both
the $\bold\Phi$- and $\bold\Xi$-function spaces.) Unlike the 
$\bold\Phi_{Qpq}$, which is not closed regarding the action of 
$\Delta_X$, the $\bold\Xi_{Qpq}$ is 
closed under the action of the
complete operator $\Delta$. Exactly this is the primarily reason why
these subspaces have been introduced here. 
Several larger $L^2$-spaces invariant under the
actions of the Laplacians can also be introduced.  
Invariant space 
$\bold\Xi_{Q,p+q=n}=\sum_{q=0}^n\bold\Xi_{Q,p=n-q,q}$,
defined for a fixed $n\in \mathbb N$,
will be thoroughly investigated later.   

Natural operator,
$
\omega_{Q\tilde Qpq}:
\bold \Phi_{Qpq}\to
\bold \Phi_{\tilde Qpq}
$
and
$
\omega_{Q\tilde Qpq}:
\bold \Xi_{Qpq}\to
\bold \Xi_{\tilde Qpq}
$, 
intertwining the corresponding function spaces defined by
two distinct unit vectors 
$Q$ and $\tilde Q$ can also be introduced. Note that this
function spaces are defined on the same H-type group. The existence
of these operators will give a clear explanation for the puzzle:
How can the spectrum ``ignore'' the isometries? 
Since spherical harmonics 
$\Pi_X(\Theta_Q^p\overline\Theta_Q^q)$ are mapped to the same type
of harmonics and 
\begin{eqnarray}
\label{MF}
\bold M\mathcal F_{Qpq}(\phi )=
\mathcal F_{Qpq}((q-p)|V|\phi )\, ,\,
\Delta_Z\mathcal F_{Qpq}(\phi )=
\mathcal F_{Qpq}(-|V|^2\phi ),\\
\label{HMF}
\bold M\mathcal{HF}_{Qpq}(\phi )=
\mathcal{HF}_{Qpq}((q-p)|V|\phi )\, ,\,
\Delta_Z\mathcal {HF}_{Qpq}(\phi )=
\mathcal {HF}_{Qpq}(-|V|^2\phi ),
\end{eqnarray}
all these operators 
obviously intertwine
the corresponding Laplacians. The second line follows from the first
one by the commutativity of operator $\bold D_V\bullet$ with the
projection $\Pi_X$. Also note that these formulas prove the intertwining
property also for $\kappa_{Qpq}$, since the same formulas are valid
also for an other group 
$H^{(a^\prime ,b^\prime )}_l$, where
$a^\prime +b^\prime=a+b$, in the family.

By these formulas, both the $\kappa_Q$- and the 
$\omega_{Q\tilde Q}$-intertwining operators can be derived from
point transformations of the form 
$(K_Q,id_Z)$ resp. $(O_{Q\tilde Q},id_Z)$,
where $id_Z$ is the identity on the Z-space and the first operators
are orthogonal transformations on the X-space. The latter ones transform
the subspace $S_Q$, spanned by $Q$ and all  $J_{Z_0}(Q)$, 
to $S^\prime_Q$ resp. $S_{\tilde Q}$. This part of the 
transformation is uniquely determined, whereas, between the complement
spaces it can be arbitrary orthogonal transformation.

This observation easily proves that the {\it Dirichlet condition} is 
intertwined
by each of the above operators on each invariant subspaces considered
above. The same statement can be proved also for the Neumann conditions.

{\bf The isospectrality} on ball- and sphere$\times$ball-type manifolds 
of nilpotent groups can be established in two different ways. 
In the first one the $L_{\bold C}^2$-Hilbert space
is decomposed into direct sum of invariant subspaces described
above, then, the $\kappa$-type intertwining operators on these subspaces
define global intertwining operator. While in the periodic case there
is only one reasonable decomposition, compelled by the 
Fourier-Weierstrass decomposition defined by the Z-lattice $\Gamma$,
in this case this splitting has to be constructed. Different
splittings define different intertwining operators.

To make the first step in this decomposition, note that the whole $L^2$ 
function space is spanned by the
$\Delta$-invariant sub-spaces
$
\bold \Xi_{Qpq}
$, considered for all $Q$, $p$, and $q$.
This statement can be established by proving that the invariant space
$\bold \Xi_{n=p+q}=\sum_{Q}\bold \Xi_{Qn=p+q}$,
defined for fixed $n=p+q$, is nothing but the space spanned by
the functions of the form $\phi(|X|)H^{(n)}(X)\varphi (Z)$, where  
$H^{(n)}$
is an $n^{th}$-order homogeneous harmonic polynomial
and $\varphi$ is an $L^2$-function. This space is spanned by
the invariant spaces   
$\bold \Xi_{pq}=\sum_{Q}\bold \Xi_{Qpq}$ defined for fixed $p$ and
$q$ satisfying $p+q=n$.
The number of independent $Q$'s  
such that the corresponding subspaces span this latter total space is 
\begin{equation}
d_{pq}=
\binom{p+k-1}{k-1}
\binom{q+k-1}{k-1},
\end{equation}
where $k$ is the dimension of the X-space. Furthermore, one can find 
$d_n=\binom{n+k-1}{k-1}$ number of independent 
vectors, $Q_1,\dots ,Q_{d_n}$, such that for any fixed $p$ and $q$
the $\bold \Xi_{pq}$ is the direct sum of subspaces 
$
\{\bold\Xi_{Q_1pq},\dots ,
\bold\Xi_{Q_{d_{pq}}pq}\}.
$
By considering all $p$ and $q$ satisfying $p+q=n$, one has a 
decomposition for $\Xi_{n=p+q}$.
One can also prove that a function 
$f\in\bold \Xi_{n=p+q}$ satisfies a given boundary condition
if and only if all of the component functions in the decomposition 
of $f$ with respect to this direct sum splitting satisfy 
the boundary condition.  
Thus, by intertwining with each
$\kappa_{Q_ipq}$, one has a global
intertwining between $\bold\Xi_{n=p+q}$ and 
$\bold\Xi^{\prime}_{n=p+q}$ which intertwines both the 
Laplacians and the boundary conditions. 

In the other isospectrality proof the intertwining operators 
$\omega_{Q\tilde Qpq}$ are applied.
The elements, $\{\lambda_{n=p+q,i}\}$, 
of the spectrum appear on each  
$
\bold \Xi_{Q,n=p+q}
$ with multiplicity, say $m_{n=p+q,i}$. Then the 
multiplicity regarding the whole $L^2$-space can be determined by
the dimensions described above. 
On the other hand, for  
$Q\in\mathbb R^{r(l)a}$, the isospectrality 
obviously follows from $\bold\Xi_{Qpq}=\bold\Xi^{\prime}_{Qpq}$. 
This proves the desired isospectrality completely.

The latter proof clearly demonstrates that the spectrum can really
be ``ignorant'' of the isometries. In fact, there is a subgroup, 
$\bold {Sp}(a)\times \bold{Sp}(b)$, 
of isometries
on a Heisenberg-type group $H^{(a,b)}_3$ which acts as the identity
on the Z-space.
Note that 
these isometries act transitively on the
X-spheres 
of $H_3^{(a+b,0)}$, 
implying the intertwining property for 
$\omega_{Q\tilde Q,n=p+q}$. Though the isometries
are not transitive on the X-spheres of the other members of 
the isospectrality family, yet, the
$\omega_{Q\tilde Q,n=p+q}$ is
still an intertwining operator on its own, 
without the help of the isometries.

{\bf Intertwining on the boundaries.} Since they intertwine the
Dirichlet boundary condition,
the above intertwining operators induce, by restrictions, maps of 
functions defined on the boundary manifolds. It turns out
that these maps are indeed intertwining operators 
on the boundary manifolds, establishing the isospectrality there too.
Then, one has new striking examples on sphere\-$\times$\-sphere-type
manifolds. In fact, the complete {\bf group of isometries} on 
sphere$\times$sphere-type manifolds of 
$H^{(a,b)}_3$ 
is the semi direct product
$(\bold {Sp}(a)\times \bold{Sp}(b))\cdot SO(3)$, 
where the action of $SO(3)$, represented by the unit quaternionic 
numbers, $q$, is defined by 
$
(X_1,\dots ,X_{a+b},Z)\to
(qX_1\overline q,\dots ,qX_{a+b}\overline q,qZ\overline q).
$
This complete group acts transitively on the 
sphere\-$\times$\-sphere-type
manifolds of
$H^{(a+b,0)}_3$, 
while the other manifolds in the isospectrality family
are locally inhomogeneous. Also the dimensions of the isometry groups
belonging to the members of the family are distinct. 

The isospectrality theorem naturally extends to 
the {\bf solvable extensions}. 
The Laplacians on the ambient- and boundary-manifolds,
furthermore, the normal vectors to the boundaries are described
in formulas (1.12), (3.30), and (3.29) of \cite{sz3}.
The generator functions
are of the form $\phi(|X|,t,V)$ in this case, i. e., the 
intertwining is led back to the nilpotent group. The group
of isometries acting on the sphere$\times$sphere-type manifolds of 
$SH^{(a,b)}_3$, where $ab\not =0$, is 
$(\bold {Sp}(a)\times \bold{Sp}(b))\cdot SO(3)$,
while on 
$SH^{a+b,0)}_3$ it is $\bold{Sp}(a+b)\cdot\bold{Sp}(1)$. As
in the nilpotent case, this yields the
same statement regarding the homogeneity property of the manifolds 
also in the solvable isospectrality family.
Thus we have
\begin{theorem}
The ball- and sphere-type manifolds, determined by the same function
$\varphi (|X|,|Z|)$ on the members of an isospectrality family  
$H^{(a,b)}_l$, are isospectral. On ball-type manifolds
this statement is proved by 
constructing, first, 
intertwining operators 
$
\kappa_{Qpq}:
\bold\Xi_{Qpq}\to
\bold\Xi_{Qpq}^{\prime}
$ 
between Laplace-invariant subspaces, then, by choosing a complete
independent subspace-system,
a global intertwining between the $L^2$-spaces is established.
In an other proof operators
$
\omega_{Q\tilde Qpq}:
\bold\Xi_{Qpq}\to
\bold\Xi_{\tilde Qpq},
$ 
associating the functions of the same H-type group to each other,
are constructed. This proof better reveals how does the spectrum 
``ignore'' the isometry group. On the sphere-type boundary manifolds
the isospectralities are established by double sided intertwining
operators.

This isospectrality statement extends
to $\sigma$-equivalent 2-step nilpotent groups defined by totally
anticommutative endomorphism spaces as well as to the ball- and 
sphere-type manifolds of the solvable extensions of these groups. 
These examples include the striking examples,
constructed on the geodesics spheres of $SH^{(a,b)}_3$, where one of 
the family-members is homogeneous while 
the others are locally inhomogeneous.

New examples are constructed on sphere$\times$ball- and
sphere$\times$sphere-type manifolds. Among them two are particularly
interesting. In fact, the isospectrality family of 
sphere$\times$sphere-type manifolds, constructed both on  
$H^{(a,b)}_3$ and
$SH^{(a,b)}_3$, the metric is homogeneous for 
the manifold belonging
to the pair $(a+b,0)$ or $(0,a+b)$, while the metrics satisfying 
$ab\not =0$ are locally inhomogeneous. Also the dimensions of
groups of isometries acting on the members are different. 
\end{theorem} 

{\bf Generalizations.} General 2-step nilpotent Lie groups are
defined by general endomorphism spaces, $J_{\bold z}$, of skew
endomorphisms acting on the X-space. A $\sigma$-deformation of the
endomorphism space is defined by an orthogonal involutive transformation
$\sigma$ on the X-space which commutes with all endomorphisms
from $J_{\bold z}$. The $\sigma$-deformed endomorphism space consists
of endomorphisms $\sigma J_Z$. Note that the family of Cliffordian 
endomorphism spaces $J^{(a,b)}_l$ defined by the same $a+b$ and $l$
consists of $\sigma$-equivalent endomorphism spaces. In 
papers \cite{sz2,sz3} 
the isospectrality is stated on the ball- and sphere-type 
domains of such 2-step nilpotent Lie groups and their solvable
extensions whose endomorphism spaces 
contain at least one anticommutator. However, the above isospectrality
constructions extend to the ball-, sphere-, 
sphere$\times$ball-, and sphere$\times$sphere-type manifolds
of two-step nilpotent Lie groups and their solvable extensions
which are defined by $\sigma$-equivalent endomorphism spaces. This
extension is non-trivial. The above ideas work out just for
H-type groups and their solvable extensions.
The extensions require substantial modifications both in the formulas
defining the intertwining operators and in investigating the boundary
conditions.

\section{Normal de Broglie Geometry}

\noindent{\bf Infinities in Quantum Electrodynamics (QED).} 
The problem of infinities (divergent integrals), which is 
present in calculations since the early days of quantum field theory
(Heisenberg-Pauli (1929-30)) or elementary particle physics
(Oppenheimer (1930) and Waller (1930) in electron theory), 
is treated by {\it renormalization} in the current theories. 
This perturbative tool provides
the desired finite quantities by differences of infinities. 
This problem is the legacy of controversial concepts such as 
{\it point mass}
and {\it point charge} of classical electron
theory, which provided the first warning that a point 
electron will have infinite electromagnetic
self-mass: the mass $e^2/6\pi ac^2$ for a surface distribution
of charge with radius $a$ blows up for $a\to 0$. In quantum
field theory the Hamiltonian of the field is proportional to this 
electromagnetic self-mass. This is why this infinity launched one 
of the deepest crisis's in the history of physics.

The infinities, related to the divergence of the summations
over all possible distributions of energy/momentum of the virtual 
particles, mostly appear in the form of infinite
traces of kernels such as the Wiener-Kac kernel $e^{-tH}$
or the Dirac-Feyn\-man kernel $e^{-tH\mathbf i}$. 
The WK-kernel provides the fundamental solution of the heat equation
while the DF-kernel provides the fundamental
solution of the Schr\"odinger equation.  
The infinite trace assigns infinite measures to 
physical objects such as
self-mass, self-charge, e.t.c.. 
Because of the divergent integrals appearing
in its construction,
also the Feynman measure, which is analogous 
to the well defined Wiener-Kac measure on the path-spaces, 
requires renormalization. 

This paper offers a new non-perturbative approach to this problem.
The main idea in this approach can be briefly described as follows.
In the first step the quantum Hilbert space 
$\mathcal H$ 
(on which the quantum
Hamilton operator $H$ is acting) is decomposed into the direct
sum of $H$-invariant subspaces, called
Zeeman zones. Then all the operator-actions, 
such as the the heat- or Feynman-flows, 
are considered on these invariant subspaces separately. 
It turns out that both the Wiener-Kac and Dirac-Feynman kernels are 
of the trace class on each zone, furthermore, both define the 
corresponding zonal measures on the path-spaces rigorously.
\medskip

\noindent{\bf Introducing the zones.}
The Hilbert space, $\mathcal H$, whose zonal decomposition 
is to be established is the space of complex valued 
$L^2$-functions defined on the X-space. 
It is isomorphic to the weighted space defined by the Gauss density 
$d\eta_\lambda (X)=e^{-\lambda |X|^2}dX$, which 
is spanned by the complex valued polynomials. 

Next the Hilbert space is interpreted in this way.
The natural {\it complex Heisenberg group 
representation} on $\mathcal H$ is defined by  
\begin{equation}
\label{rho}
\rho_{\mathbf c} (z_i)(\psi )=(-\partial_{\overline{z}_i}+
\lambda z_i\cdot )\psi\quad ,\quad
\rho_{\mathbf c} (\overline{z}_i)(\psi )=\partial_{z_i}\psi , 
\end{equation}
where $\{z_i\}$ is a complex coordinate system on the X-space.
This representation is reducible. In fact, it is irreducible on the
space of holomorphic functions, where it is called Fock representation.
Besides the holomorphic subspace there are infinitely many
other irreducible invariant subspaces. In the literature
only the Fock representation, defined on the space of 
holomorphic functions, is well known. The above representation is
called {\it extended Fock representation}. In the function operator
correspondence, this representation associates operator (\ref{Zee_int})
to the Hamilton function of an electron orbiting in constant 
magnetic field.

The zones are defined in two different ways. First, they can be defined
by the invariant subspaces of representation (\ref{rho}). The actual 
construction uses Gram-Schmidt orthogonalization. On the complex plane
$\mathbf v=\mathbf C$, which corresponds to the 2D-Zeeman operator 
(\ref{Zee_int}),
the first zone, $\mathcal H^{(0)}$, is the holomorphic zone spanned
by the holomorphic polynomials. To construct the second zone
one considers the function space $G^{(1)}$ consisting of functions 
of the form $\overline zh$, where the $h$ is a holomorphic polynomial.
Then $\mathcal H^{(1)}$ is the orthogonal component of $G^{(1)}$
to the previous holomorphic zone. E. t. c., one has
all zones, $\mathcal H_\lambda^{(a)}$, by continuing with the
Gram-Schmidt orthogonalization applied to function spaces $G^{(a)}$,
which are spanned by functions of the form $\overline z^ah$. 
The zone index $a$ indicates the maximal number
of the antiholomorphic coordinates $\overline z$ in the polynomials
spanning the zone. 

In the 2D-case all the zones are irreducible under the action of the 
extended Fock representation. In the higher dimensions, however,
the Gram-Schmidt process results reducible zones, called 
{\it gross zones.} More precisely, the holomorphic zone is
always irreducible and the gross zones of higher indexes 
decompose into {\it irreducible zones} 
which can also be explicitly described.

The second technique defines the very same zones by  
computing the spectrum and the corresponding 
eigenfunctions explicitly.  According to these computations, the
eigenfunctions appear in the form
\begin{equation}
h^{(p,\upsilon)}(X)=
H^{(p,\upsilon)}(X) 
e^{-\lambda |X|^2/2}
\end{equation}
with the corresponding eigenvalues
\begin{equation}
\label{eigv1}
-((4p+k)\lambda +4k\lambda^2),  
\end{equation}
where $p$ resp. $\upsilon$ are the holomorphic resp. antiholomorphic
degrees of polynomial $H^{(p,\upsilon )}$. 
Numbers $l=p+\upsilon$ and $m=2p-l$ are
called azimuthal and magnetic quantum numbers respectively. The
above function is an eigenfunction also of the magnetic dipole
moment operator with eigenvalue $m$.

In this new way the zones are created such that the above 
eigenfunction falls into the zone with index $\upsilon$. 
According to the formula $\upsilon ={1\over 2}(l-m)$, 
the zones are determined by the magnetic quantum number $m$.
Thus a zone exhibits the magnetic state of the particle.

Note that eigenvalues (\ref{eigv1}) are independent of the 
antiholomorphic index and they depend just on the holomorphic index. 
As a result, each eigenvalue has infinite multiplicity. 
On the irreducible zones, however,
each multiplicity is $k/2$. (Here we suppose that there is only one 
parameter $\lambda$ involved, meaning that the particles are identical.
If the particles (i. e., the $\lambda_i$'s) 
are properly distinct, then the 
multiplicity is $1$ on each zone.) Moreover, two irreducible zones 
are isospectral.
\medskip

\noindent{\bf Introducing the point-spreads by projection kernels.} 
It is remarkable that all the 
important objects such as the projections onto the zones,
the zonal heat- and Feynman-kernels and their well defined trace,
the zonal partition functions (defined 
with no regularization), and the several Feynman-Kac type formulas 
can explicitly be computed \cite{sz5,sz6}. 

First the {\it projection operators}, $\delta^{(a)}$, 
onto the zones $\mathcal H^{(a)}$ are established. 
If $\{\varphi_i^{(a)}\}_{i=1}^\infty$ is an 
orthonormal basis in $\mathcal H^{(a)}$, 
then the corresponding projection can be
formally defined as convolution with the kernel 
\begin{equation}
\delta^{(a)}(z,w)=\sum_i \varphi^{(a)}_i(z)
\overline{\varphi^{(a)}_i(w)},
\end{equation}
where $z$ and $w$ represent complex vectors on 
$\mathbf C^{k\over 2}=\mathbf R^k$.
Interestingly enough, these self-adjoint operators are integral
operators having a smooth Hermitian integral kernels.
These kernels can be interpreted as restrictions of the 
global Dirac delta distribution, 
$\delta_z(w)=\sum \varphi_i(z)\overline{\varphi}_i(w)$, onto the
zones. They have the following
explicit form
\begin{equation}
\delta_{\lambda z}^{(a)}(w)
= {\lambda^{k/2}\over \pi^{k/2}}L^{((k/2)-1)}_a(\lambda |z-w|^2)
e^{\lambda (z\cdot\overline{w}-{1\over 2}(|z|^2+ |w|^2)},
\end{equation}
where $L_a^{((k/2)-1)}(t)$ is the corresponding Laguerre polynomial.
Among these kernels only the first one, the projection kernel
onto the holomorphic zone, is well known. It is nothing but the
Bergman kernel. The new mathematical feature of these formulas is
that they are explicitly determined regarding each zone and not just
for the holomorphic zone.

These kernels represent one of the most important concepts
in this theory. They can be interpreted such that, 
on a zone, a point particle appears as a spread
described by the above wave-kernel. 
Note that how these kernels, called zonal point-spreads, are derived
from the one defined for the holomorphic zone. This holomorphic 
spread is just multiplied by the radial 
Laguerre polynomial corresponding to the zone. These point-spreads
show the most definite similarity to the de Broglie waves
packets (cf. \cite{bo}, pages 61), suggesting
that a point particle concentrated
at a point $Z$ appears on a zone as an object spread around 
$Z$ as a wave-package. The wave-function is 
described by the above kernel explicitly. 

The wave-package interpretation of physical objects started out with the
de Broglie theory. This concept was finalized in the Schr\"odinger
equation. The mathematical formalism did not follow this development,
however, and the Schr\"odinger theory is built up on a mathematical
background not excluding the existence of the controversial
point objects. On the contrary, 
an electron must be considered as a point-object in
the Schr\"odinger theory as well 
(cf. Weisskopf's argument on this problem in \cite{schw, sz5}). 
An other demonstration for 
the presence of
point particles in classical theory is the duality principle, 
stating that objects manifest themselves sometimes as waves and 
sometimes
as point particles. The bridge between the two visualizations is
built up in Born's probabilistic theory,
where the probability for that that a particle, 
attached to a wave $\xi$, can be found on a domain $D$ is measured by
$\int_D\xi\overline\xi$.

These controversial point-objects, by having infinite self-mass or
self-charge attributed to them by the Schr\"odinger equation, 
launched one of the deepest crisis's in the history of physics.
In the zonal theory de Broglie's idea is established on a mathematical
level. Although the points are ostracized from this theory, 
the point-spreads 
still bear some reminiscence of the point-particles. For instance, they
are the most compressed wave-packages and all the 
other wave-functions in the 
zone can be expressed as a unique superposition of the point-spreads.
If $\xi$ is a zone-function, the above integral 
measures the probability that
the center of a point-spread is on the domain $D$. This 
interpretation restores, in some extend, the duality principle in the
zonal theory.

Function 
$
\delta^{(a)}_{\lambda Z}
\overline{\delta}^{(a)}_{\lambda Z}
$
is called the density of the spread around $Z$. By this reason, 
function 
$
\delta^{(a)}_{\lambda Z}
$
is called spread-amplitude. Both the spread-amplitude and 
spread-density generate well defined measures on the path-space.
This space consists of continuous curves 
connecting two arbitrary points.
Both measures can be constructed by the method applied in constructing
the Wiener measure.
 
The point-spread 
concept bears some remote reminiscence of Heisenberg's 
suggestion (1938) for the existence of 
a fundamental length $L$, analogously
to $h$, such that field theory was valid only for distances larger than
$L$ and so divergent integrals would be cut off at that distance.
This idea has never became an effective theory, however. 
Other distant relatives of the point-spread-concept are the 
smeared operators, i. e. those suitably averaged over small 
regions of space-time, considered by Bohr and Rosenfeld in quantum 
field theory. There are also other theories where an electron is 
considered to be extended. Most of them fail on lacking the
explanation for the question: 
Why does an extended electron not blow up? The zonal theory is checked
against this problem in \cite{sz5}, section 
(F) ``Linking to the blackbody radiation; Solid zonal particles''.
\medskip

\noindent{\bf Global Wiener-Kac and Feynman flows.}
Both definitions imply that the zones are invariant 
under the action of the Hamilton operator,
therefore, kernels such as the {\it heat (Wiener-Kac) and 
Feynman-Dirac kernels} can be restricted onto them. The 
zonal kernels are defined by these restrictions. 

Since the spectrum is discrete, also the global kernels, 
defined for the total space $\mathcal H$, can be introduced 
by the trace formula using an orthonormal basis
consisting of eigenfunctions on the whole space $\mathcal H$.
Despite of the infinite multiplicities on the global setting,
both global kernels are well defined smooth functions. 
If the Zeeman operator $H_Z$ is non-degenerated such that the
distinct non-zero parameters $\{\lambda_i\}$,  $i=1,\dots ,r$,
are defined on $k_i$-dimensional subspaces, then for the Wiener-Kac 
kernel we have
\begin{eqnarray}
\label{d_1gamm}
d_{1\gamma} (t,X,Y)=e^{-tH_Z}(t,X,Y)=\\
=\prod \big({\lambda_i\over 2\pi sinh(\lambda_i t)}\big)^{k_i/2}e^{
-\sum{\lambda_i}({1\over 2}coth(\lambda_i t)
|X_i-Y_i|^2+\mathbf i
\langle X_i,J (Y_i)\rangle}.\nonumber
\end{eqnarray}
This kernel satisfies the Chapman-Kolmogorov identity as well as the
limit property $\lim_{t\to 0_+}d_1(t,X,Y)=\delta(X,Y)$, however, 
it is not of the trace class. Thus functions such as 
the partition function or the zeta function are not defined in the
standard way. Note that by regularization (renormalization) only
well defined relative(!) partition and zeta functions are introduced.

The explicit form of the global Feynman-Dirac kernel is
\begin{eqnarray}
d_{\mathbf i}(t,X,Y)=e^{-t\mathbf iH_Z}(t,X,Y)=\\
=\prod\big({\lambda_i\over 2\pi
\mathbf i sin(\lambda_i t)}\big)^{k_i/2}e^
{\mathbf i\sum{\lambda_i}\{
{1\over 2} cot(\lambda_i t)|X_i-Y_i|^2-\langle X_i,J(Y_i)\rangle \}}.
\nonumber
\end{eqnarray}
Since for fixed $t$ and $X$ function $d_{\mathbf i}(t,X,Y)$ is neither
$L^1$- nor $L^2$-function of the variable $Y$, 
the integral required
for the Chapman-Kolmogo\-rov identity is not defined for this kernel. 
Neither is it of the 
trace class. Nevertheless, it satisfies the 
above limit property. Thus the
constructions with the Feynman kernel lead to divergent 
integrals in the very first step.

It is well known in the history that 
Kac, who tried to understand Feynman, was able to introduce a 
rigorously defined
measure on the path-spaces only by the kernel $e^{-tH}$. 
This measure was, actually, established earlier by Wiener for the 
Euclidean Laplacian $\Delta_X$. 
Note that the heat kernel involves a Gauss density which makes this
constructions possible. Whereas, the Feynman-Dirac 
kernel does not involve
such term. This is why no well defined constructions can be
carried out with this kernel. One can strait out all this
difficulties, however, by considering these 
constructions on the zones separately.
\medskip

\noindent{\bf Zonal Wiener-Kac and Feynman flows.}
Also the zonal WK- resp. FD-kernels are well defined smooth functions.
The gross {\bf zonal
Wiener-Kac kernels} are of the trace class, which can be
described, along with their {\bf partition functions}, by the following
explicit formulas. 
\begin{eqnarray}
\label{d_1^a}
d_{1}^{(0)}(t,X,Y)=
\prod\big({\lambda_i e^{-\lambda_i t}\over \pi}\big)^{k_i\over 2}
e^{\sum \lambda_i(-{1\over 2} (|X_i|^2+|Y_i|^2)+ e^{-2\lambda_i t}
\langle X_i,Y_i+\mathbf iJ(Y_i)\rangle )},
\\
d_{1}^{(a)}(t,X,Y)=
(L^{({k\over 2}-1)}_a(\sum\lambda_i |X_i-Y_i|^2)+LT_1^{(a)}
(t,X,Y))d_1^{(0)}(t,X,Y),
\end{eqnarray}
where 
$LT_1^{(1)}$ 
is of the form 
\begin{eqnarray}
LT_1^{(1)}(t,X,Y))=
(1-e^{-2t})
lt_1^{(1)}(t,X,Y))=
\\
(1-e^{-2t})
(|X|^2+|Y|^2-1-
(1+e^{-2t})\langle X,Y+\mathbf iJ(Y)\rangle ) \nonumber
\end{eqnarray}
and for the general terms,
$LT_1^{(a)}$,
recursion formula can be established. Furthermore,
\begin{eqnarray}
\mathcal Z_1^{(a)}(t)=
Trd_{1}^{(a)}(t)
={a+(k/2)-1\choose a}\prod {e^{-{k_i\lambda_i t\over 2}}\over
(1-e^{-2\lambda_i t})^{k_i\over 2}}=
TrD_{1}^{(a)}(t),
\end{eqnarray}
where 
$
D_{1}^{(a)}(t,X,Y)=
L^{({k\over 2}-1)}_a(\sum\lambda_i|X_i-Y_i|^2)
d_{1}^{(0)}(t,X,Y)
$
is the dominant zonal kernel. The remaining long term kernel
in the WK-kernel vanishes for 
$\lim_{t\to 0_+}$ and is of the 0 trace
class.
The zonal WK-kernels satisfy the Chapman-Kolmogorov identity 
along with the limit property
$\lim_{t\to 0_+}d_1^{(a)}=\delta^{(a)}$. 

Similar statements can be established regarding the zonal
DF-flow. The gross {\bf zonal 
Dirac-Feynman kernels} are of the trace class which, together with their
{\bf partition functions}, can be described 
by the following explicit formulas.

\begin{eqnarray}
\label{d_i^a}
d_{\mathbf i}^{(0)}(t,X,Y)=
\prod\big({\lambda_i e^{-\lambda_i t\mathbf i}\over \pi}
\big)^{k_i\over 2}
e^{\sum \lambda_i(-{1\over 2} (|X_i|^2+|Y_i|^2)+ 
e^{-2\lambda_i t\mathbf i}
\langle X_i,Y_i+\mathbf iJ(Y_i)\rangle )},
\\
d_{\mathbf i}^{(a)}(t,X,Y)=
(L^{({k\over 2}-1)}_a(\sum\lambda_i |X_i-Y_i|^2)+LT_{\mathbf i}^{(1)}
(t,X,Y))d_{\mathbf i}^{(0)}(t,X,Y),
\end{eqnarray}
where 
$LT_{\mathbf i}^{(1)}$ 
is described by 
\begin{eqnarray}
LT_{\mathbf i}^{(1)}(t,X,Y))=
(1-e^{-2t\mathbf i})
lt_{\mathbf i}^{(1)}(t,X,Y))=
\\
(1-e^{-2t\mathbf i})
(|X|^2+|Y|^2-1-
(1+e^{-2t\mathbf i})\langle X,Y+\mathbf iJ(Y)\rangle ).\nonumber
\end{eqnarray}
and a
general long term, 
$LT_{\mathbf i}^{(a)}$,
can be defined recursively. Furthermore, 
\begin{eqnarray}
\mathcal Z_{\mathbf i}^{(a)}(t)=
Trd_{\mathbf i}^{(a)}(t)
={a+(k/2)-1\choose a}\prod {e^{-{k_i\lambda_i t\mathbf i\over 2}}\over
(1-e^{-2\lambda_i t\mathbf i})^{k_i\over 2}}=
TrD_{\mathbf i}^{(a)}(t),
\end{eqnarray}
where 
$
D_{\mathbf i}^{(a)}(t,X,Y)=
L^{({k\over 2}-1)}_a(\sum\lambda_i|X_i-Y_i|^2)
d_{\mathbf i}^{(0)}(t,X,Y)
$
is the dominant kernel. The remaining longterm term in the zonal
DF-kernel is of the 0 trace class.

The zonal DF-kernels are zonal fundamental solutions 
of the Schr\"o\-dinger equation 
$(\partial_t+\mathbf i(H_{Z})_{X})d^{(a)}_{\mathbf i\gamma}(t,X,Y)=0$,
satisfying the Chapman-Kolmogorov identity as well as 
the limit property $\lim_{t\to 0_+}d_{\mathbf i}^{(a)}=\delta^{(a)}$.
 
On the zones the WK- and the FD-kernels are not just of the 
trace class. They both define, rigorously, 
complex zonal measures, 
the {\bf zonal Wiener-Kac measure} 
$dw_{1xy}^{T(a)}(\omega)$ 
and the {\bf zonal Feynman measure}
$dw_{\mathbf ixy}^{T(a)}(\omega)$,
on the space of continuous curves $\omega :[0,T]\to \mathbb R^k$ 
connecting two points $x$ and $y$.
The existence of zonal WK-measure
is not surprising, since this measure exists even for the 
global setting.
However, the trace class property is a new feature, indeed. In
case of the zonal Feynman measure both the trace class property 
and the existence of rigorous zonal Feynman measures are new features. 
Note that also the zonal DF-kernels
involve a Gauss kernel which makes these constructions well defined.

\end{document}